\documentclass[12pt,reqno]{amsart}
\usepackage{amssymb}
\usepackage{graphicx}
\usepackage{amsmath}
\usepackage{a4wide}

\newtheorem{theorem}{Theorem}

\newtheorem{corollary}[theorem]{Corollary}
\newtheorem{lemma}[theorem]{Lemma}
\newtheorem{proposition}[theorem]{Proposition}

\theoremstyle{remark}
\newtheorem{remark}[theorem]{Remark}

\def\leq{\leqslant}
\def\geq{\geqslant}
\def\al{\alpha}
\def\be{\beta}
\def\ga{\gamma}
\def\C{\mathbb C}
\def\F{\mathbb F}

\def\R{\mathbb R}
\def\Q{\mathbb Q}
\def\Z{\mathbb Z}
\def\dmax{d_{\max}}

\def\ep{\varepsilon}

\def\cc{,\dots,}
\def\l{\ell}
\def\om{\omega}
\def\Ql{\Q(\om_\l)}
\def\MW{Mordell-\mbox{\kern-.16em}Weil}

\newcommand{\Gal}{{\operatorname{Gal}}}
\newcommand{\GL}{{\operatorname{GL}}}
\newcommand{\SL}{{\operatorname{SL}}}
\newcommand{\ST}{{\operatorname{ST}}}
\newcommand{\Qbar}{{\overline{\Q}}}

\newcommand{\Fbar}{{\overline{\F}}}
\newcommand{\kbar}{{\overline{k}}}
\newcommand{\Stab}{{\operatorname{Stab}}}
\newcommand{\tensor}{\otimes}
\newcommand{\mm}{{\mathfrak m}}

\begin{document}
\title[Dimension of algebraic numbers]{The conjugate dimension 
					of algebraic numbers}
\author{Neil Berry} 
\address{School of Mathematics \\
University of Edinburgh\\
James Clerk Maxwell Building\\
King's Buildings, Mayfield Road\\
Edinburgh EH9 3JZ\\
Scotland, U.K.}
\email{ neilb@maths.ed.ac.uk}
\author{Art\= uras Dubickas}
\address{Department of Mathematics and Informatics\\ Vilnius University\\
Naugarduko 24, Vilnius 03225\\ Lithuania}
\email{arturas.dubickas@maf.vu.lt}
\author{ Noam D. Elkies}
\address{Department of Mathematics\\ Harvard University\\ 
			Cambridge, MA 02138\\ USA}
\email{elkies@math.harvard.edu}
\author{ Bjorn Poonen}
\address{Department of Mathematics\\ University of California \\ 
			Berkeley, CA 94720-3840\\ USA}
\email{poonen@math.berkeley.edu}
\author{ Chris~Smyth}
\address{School of Mathematics \\
University of Edinburgh\\
James Clerk Maxwell Building\\
King's Buildings, Mayfield Road\\
Edinburgh EH9 3JZ\\
Scotland, U.K.}
\email{c.smyth@ed.ac.uk}
\date{4 May, 2004}

\begin{abstract}
We find sharp upper and lower bounds for the degree of an algebraic number in
terms of the $\Q$-dimension of the space spanned by its conjugates. 
For all but seven nonnegative integers $n$ the
largest degree of an algebraic number whose conjugates
span a vector space of dimension $n$ is equal to $2^n n!$. 
The proof, which covers also the seven exceptional cases,
uses a result of Feit on the maximal order of finite subgroups of $\GL_n(\Q)$;
this result depends on the classification of finite simple groups.
In particular, we construct an algebraic number of degree $1152$ 
whose conjugates span a vector space of dimension only $4$.

We extend our results in two directions. We consider the problem when
$\Q$ is replaced by an arbitrary field, and prove some general results.
In particular, we again obtain sharp bounds when the ground field
is a finite field, or a cyclotomic extension $\Ql$ of $\Q$.
Also, we look at a multiplicative version of the problem
by considering the analogous rank problem for the multiplicative group
generated by the conjugates of an algebraic number.
\end{abstract}

\maketitle

\section{Introduction}\label{}  
Let $\Qbar$ be an algebraic closure of 
the field $\Q$ of rational numbers, 
and let $\al \in \Qbar$.
Let $\al_1,\dots,\al_d \in \Qbar$ 
be the conjugates of $\al$ over $\Q$, with $\al_1=\al$.
Then $d$ is the degree $d(\al):=[\Q(\al):\Q]$,
the dimension of the $\Q$-vector space spanned by the powers of $\al$.
In contrast, we define the {\em conjugate dimension} $n=n(\al)$ of $\al$
as the dimension of the $\Q$-vector space spanned 
by $\{\al_1,\dots,\al_d\}$.

In this paper we compare $d(\al)$ and $n(\al)$.
By linear algebra, $n \leq d$.
If $\al$ has nonzero trace and has Galois group
equal to the full symmetric group $S_d$,
then $n=d$ (see~\cite[Lemma~1]{smyth1986}).
On the other hand, it is shown in \cite{dubickas2003}  that $n$ can be 
as small as $\lfloor \log_2 d \rfloor$. 
It turns out that $n$ can be even smaller. 
Our first main result gives the minimum and maximum values of $d$
for fixed $n$.

\begin{theorem}\label{T:main}
Fix an integer $n \geq 0$.
If $\al \in \Qbar$ has $n(\al)=n$,
then the degree $d=d(\al)$ satisfies
$n \leq  d \leq \dmax(n)$, 
where $\dmax(n)$ is defined by Table~\ref{Ta1},
equalling $2^n n!$ for all $n \notin \{2,4,6,7,8,9,10\}$.
Furthermore, for each $n \geq 1$,
there exist $\al \in \Qbar$ attaining the lower and upper bounds.
\end{theorem}

We refer to the $n$ with $\dmax(n)\ne 2^nn!$ as  {\em exceptional}.
To attain $d=\dmax(n)$, we will use $\al$ 
for which the extension $\Q(\al)/\Q$ is Galois with
 Galois group isomorphic to a maximal-order finite subgroup $G$ of $\GL_n(\Q)$ 
given in Table~\ref{Ta1}. 
 
\begin{table}[ht]
\begin{center}
\bigskip\begin{tabular} { r | r  | l | r}

$ n$  & $\dmax(n) / (2^n n!) $ & Maximal-order subgroup $G$
  & $\dmax(n)=\# G$\\
     \hline
     2  & 3/2    & $W(G_2)$                   & 12\\
     4  & 3      & $W(F_4)$                   & 1152 \\
     6  & 9/4    & $\langle W(E_6),-I\rangle$ & 103680\\
     7  & 9/2    & $W(E_7)$                   & 2903040\\
     8  & 135/2  & $W(E_8)$                   & 696729600 \\
     9  & 15/2   & $W(E_8)\times W(A_1)$      & 1393459200\\
    10  & 9/4    & $W(E_8)\times W(G_2)$      & 8360755200\\ 
all other $n$ & 1 & $W(B_n) = W(C_n) = (\Z/2\Z)^n \rtimes S_n$ & $2^n n!$ \\
    \end{tabular}
    \vspace*{1ex}
    \caption{Maximal-order finite subgroups of $\GL_n(\Q)$} \label{Ta1}
\end{center}\end{table}

The groups $W(\cdot)$ are the Weyl groups of classical Lie algebras 
acting on their maximal tori (see for instance~\cite{humphreys1990}). 
They are all reflection groups: 
each is generated by those elements that
act on~$\Q^n$ by reflection in some hyperplane.
For the standard fact that the negative identity
matrix $-I$ is not in $W(E_6)$, see for instance \cite[p.~82]{humphreys1990}.
In particular, $W(B_n)=W(C_n)=(\Z/2\Z)^n \rtimes S_n$
is better known as the {\em signed permutation group},
the group of $n \times n$ matrices with entries in $\{-1,0,1\}$
having exactly one nonzero entry in each row and each column.

Feit~\cite{feit1996} proved that for each~$n$ a subgroup of $\GL_n(\Q)$ 
of maximal finite order is conjugate
to the group given in Table~\ref{Ta1}.
(The paper~\cite{feit1996} is just a statement of results --- no proofs.)
Feit's result uses unpublished work of Weisfeiler 
depending on the classification theorem for 
finite simple groups (see also \cite[p.~185]{kuzmanovich-pavlichenkov2002}). 
See
\begin{center}
  {\tt http://weisfeiler.com/boris/philinq-8-28-2000.html}
\end{center} 
for the sad tale of Weisfeiler's disappearance.

The inequality $d \leq \dmax(n)$ comes from studying the
span of $\{\al_1,\dots,\al_d\}$
as a representation of $\Gal(\Q(\al_1,\dots,\al_d)/\Q)$.
To prove the existence of examples where this upper bound is attained, 
we
\begin{enumerate}
\item 
   observe that if $G$ is one of the maximal-order finite subgroups
   of $\GL_n(\Q)$ listed in Table~\ref{Ta1}, then the $G$-invariant
   subfield $\Q(x_1,\dots,x_n)^G$ of $\Q(x_1,\dots,x_n)$
   is purely transcendental, say $\Q(f_1,\dots,f_n)$
   (whence $\Q(x_1,\dots,x_n)/\Q(f_1,\dots,f_n)$
   is a Galois extension with Galois group $G$\/),
\item apply Hilbert irreducibility to obtain a Galois extension~$K$\/ of~$\Q$
	with Galois group $G$, and
\item choose $\al \in K$\/ generating a suitable subrepresentation of $G$.
\end{enumerate}
Moreover, we give explicit examples for all $n$ 
except $6$, $7$, $8$, $9$, $10$,
and outline an explicit construction in these remaining five cases.

Many of the arguments work over base fields other than $\Q$,
so we generalize as appropriate (Theorem \ref{T:otherfields}).
In particular, Theorem~\ref{T:cyc} generalizes Theorem \ref{T:main}
by giving the minimal and maximal degrees
over any cyclotomic base field~$\Ql$.
The answers change drastically for base fields of positive characteristic: 
for instance from Theorem~\ref{T:otherfields}(v)
there are elements of a separable closure of $\F_q(t)$
of conjugate dimension~$2$ that generate Galois extensions of $\F_q(t)$
of arbitrarily large degree.
We also give in Section~\ref{multrank} some results on analogous questions
concerning the rank of the multiplicative subgroup of~$\Qbar^*$
generated by $\al_1,\ldots,\al_d$,
and its generalization over a Hilbertian field.

\section{Degree and conjugate dimension over fields in general}

\subsection{Representations}
Let $k$ be a field, and let $k^s$ be a separable closure of $k$.
If $\al \in k^s$, 
then let $d=d(\al)$ be the degree $[k(\al):k]$,
and let $n=n(\al)$ be the {\em conjugate dimension} of $\al$ (over~$k$),
defined as the dimension of the $k$-vector space $V(\al)$
spanned by the conjugates $\al_1,\dots,\al_d$ of $\al$ in~$k^s$.

\begin{proposition}\label{P:faithful} 
With notation as above,
let $K = k(\al_1,\dots,\al_d)$ and let $G = \Gal(K/k)$.
Then there exists a faithful $n$-dimensional $k$-representation of $G$.
\end{proposition}

\begin{proof}
Since $\{\al_1,\dots,\al_d\}$ is $G$-stable,
the \hbox{$G$\/-action} on $K$\/
restricts to a \hbox{$G$\/-action} on~$V(\al)$.
If $g\in G$\/ acts trivially on~$V(\al)$,
then $g$ fixes each $\al_i$,
so $g$ is the identity on $K$.
Thus $V(\al)$ is a faithful $k$-representation of $G$.
Finally, $\dim_k V(\al) = n$, by definition.
\end{proof}

A partial converse will be given in Proposition~\ref{P:subregular}
below, whose proof relies on the following representation-theoretic result.

\begin{lemma}
\label{L:faithful subrep}
Let $k$ be a field of characteristic~$0$, and let $G$ be a finite group.
Let $V$ be a $kG$-submodule of the regular representation $kG$.
Assume that $G$ acts faithfully on~$V$.
Then $V = (kG)\al$ for some $\al \in V$ with $\Stab_G(\al)=\{1\}$.
\end{lemma}

\begin{proof}
Since $k$ has characteristic~zero,
$V$\/ is a direct summand (and hence a quotient) of the regular representation, 
so the $kG$-module $V$\/ can be generated by one element.
An element $\al \in V$\/ fails to generate $V$\/ as a $kG$-module 
if and only if $\{g\al: g\in G\}$ fails to span $V$,
and this condition can be expressed in terms of the vanishing
of certain minors in the coordinates of $\al$ 
with respect to a basis of~$V$.
Thus the set $Z:=\{\,\al \in V : (kG)\al \ne V \,\}$ of such elements
is contained in the zeros of some nonzero polynomial in the coordinates.
Also, for each $g\in G - \{1\}$, the set $V^g:=\{v \in V: gv=v\}$ is
a proper subspace of $V$, since $V$\/ is faithful.
Since $k$ is infinite,
we can choose $\al \in V$ outside $Z$\/ and each $V^g$ for $g \ne 1$.
\end{proof}

\begin{remark}
We may also allow $k$\/ to have characteristic $p>0$,
as long as $p$ does not divide $\#G$ and $k$\/ is infinite.
Then $V$\/ is still a direct summand and a quotient of~$kG$,
and the same proof applies.
The hypothesis that $k$\/ is infinite cannot be removed, however,
as the following counterexample shows.  Let $k$\/ be a finite
field of characteristic~$p$, let $k'/k$\/ be a finite extension,
and take $V=k'$.  For any subgroup $G_1$ of $\Gal(k'/k)$,
let $G$\/ be the semidirect product ${k'}^* \rtimes G_1$,
which acts \hbox{$k$-linearly} on~$V$.  Then every nonzero
$\al\in V$\/ has stabilizer isomorphic to~$G_1$.
If moreover $\#G_1$ equals neither~$1$ nor a multiple of~$p$,
then $p$ does not divide $\#G$, and thus $V$\/ is a submodule of~$kG$\/
since  $V$\/ is multiplicity-free over $\kbar$\/;
but the conclusion of Lemma~\ref{L:faithful subrep} is false
because no $\al\in V$\/ has trivial stabilizer.
\end{remark}

\begin{proposition}
\label{P:subregular}
Let $k$ be a field of characteristic~$0$, and let $G$ be a finite group.
Suppose that $G=\Gal(K/k)$  for some Galois extension $K$ of $k$,
and that there is a faithful $n$-dimensional subrepresentation $V$
of the regular representation of $G$ over $k$.
Then there exists $\al \in K$ with $n(\al)=n$ and $d(\al)=[K:k]=\#G$.
\end{proposition}

\begin{proof}
By the Normal Basis Theorem, $K$, as a representation of $G$ over $k$,
is isomorphic to the regular representation.
Hence we may identify $V$\/ with a subrepresentation of $K$.
Lemma~\ref{L:faithful subrep} gives an element $\al \in V$
whose $G$-orbit has size $\#G$ and spans the $n$-dimensional space $V$.
\end{proof}

\subsection{Invariant subfields}
\label{subsec:inv}

\begin{proposition}
\label{P:invariants}
Let $G$ be one of the groups in Table~\ref{Ta1},
viewed as a subgroup of $\GL_n(\Q)$.
Then for any field $k$ of characteristic~$0$,
the invariant subfield $k(x_1,\dots,x_n)^G$ 
is purely transcendental over $k$.
\end{proposition}

\begin{proof}
We may assume $k=\Q$.
Chevalley~\cite{chevalley1955} proved that if $G$ is a finite reflection group,
then $\Q[x_1,\dots,x_n]^G = \Q[f_1,\dots,f_n]$
for some homogeneous polynomials $f_i$.
In this case, we have $\Q(x_1,\dots,x_n)^G=\Q(f_1,\dots,f_n)$ as desired.

The only remaining case is $n=6$ and $G=\langle W(E_6),-I \rangle$.
Here $\Q(x_1,\dots,x_6)^{W(E_6)}=\Q(I_2,I_5,I_6,I_8,I_9,I_{12})$
where each $I_j$ is a homogeneous polynomial of degree $j$, 
given explicitly for instance in~\cite{frame1951} (see also~\cite[p.~59]{humphreys1990}).
Moreover $-I \in G$ acts on this subfield by $I_j \mapsto (-1)^j I_j$,
so $\Q(x_1,\dots,x_6)^G = \Q(I_2,I_6,I_8,I_{12},I_5^2,I_5I_9)$.
\end{proof}

\begin{remark}
Let $G$ be a finite subgroup of $\GL_n(\R)$.
Coxeter showed~\cite{coxeter1951}
that $\R[x_1,\dots,x_n]^G$ is a polynomial ring over $\R$
in $n$ algebraically independent generators
if $G$\/ is a finite reflection group.
Shephard and Todd proved that this sufficient condition on~$G$\/
is also necessary (\cite[Thm.~5.1]{shephard-todd1954}, see also \cite[p.~65]{humphreys1990}).
For example, $G=\langle W(E_6),-I \rangle$ is not a finite reflection group,
and the $\R$-algebra 
$\R[x_1,\dots,x_6]^G = \R[I_2,I_6,I_8,I_{12},I_5^2,I_5I_9,I_9^2]$
cannot be generated by $6$ polynomials.
\end{remark}

\subsection{Hilbert irreducibility}

It is well known that the field $\Q$ is Hilbertian --- see for instance 
 \cite[Theorem~3.4.1]{serretopics}
(a form of the Hilbert Irreducibility Theorem).
This implies that Galois extensions of purely transcendental 
extensions $\Q(f_1,\dots,f_n)$ can be specialized
to Galois extensions of $\Q$ having the same Galois group 
\cite[Corollary~3.3.2]{serretopics}.

\begin{proposition}
\label{P:hilbert}
Let $k$ be a Hilbertian field.
Let a finite subgroup $G$ of $\GL_n(k)$
act on $k(x_1,\dots,x_n)$ so that
the action on the span of the indeterminates $x_i$
corresponds to the inclusion of $G$ in $\GL_n(k)$.
If the invariant subfield $k(x_1,\dots,x_n)^G$
is purely transcendental over $k$,
then there exists a finite Galois extension $K$ of $k$ with Galois group $G$.
\end{proposition}

\begin{proof}
By assumption $k(x_1,\dots,x_n)^G=k(f_1,\dots,f_n)$ for some
algebraically independent~$f_i$.
By Galois theory, $k(x_1,\dots,x_n)$ 
is a Galois extension of $k(f_1,\dots,f_n)$ with Galois group~$G$.
Now use the assumption that $k$ is Hilbertian to specialize.
\end{proof}

\begin{corollary}
\label{C:inverseGalois}
If $k$ is a Hilbertian field,
and $G$ is one of the groups in Table~\ref{Ta1},
then $G$ is realizable as a Galois group over $k$.
\end{corollary}

\begin{proof}
Combine Propositions \ref{P:invariants} and~\ref{P:hilbert}.
\end{proof}

For background material on Hilbert irreducibility 
see~\cite{schinzel2000} or~\cite{serretopics}.

\section{Degree and conjugate dimension over $\Q$}

\subsection{Proof of Theorem~\ref{T:main}}

\begin{proof}
The inequality $n \leq d$ is immediate.
Examples with equality exist
by Proposition~\ref{P:subregular} applied to the standard permutation
representation $S_n \hookrightarrow \GL_n(\Q)$,
since $S_n$ is realizable as a Galois group over $\Q$
(see~\cite[p.~42]{serretopics}, for example).

On the other hand, $d \leq \#G \leq \dmax(n)$,
where $G$ is the Galois group of $\al$ over~$k$,
because of Proposition~\ref{P:faithful},
since $\dmax(n)$ is the size of the largest finite subgroup of $\GL_n(\Q)$.

Finally, we prove that $d=\dmax(n)$ is possible for each $n \geq 1$.
Let $G$ be a maximal finite subgroup of $\GL_n(\Q)$,
as in Table~\ref{Ta1}.
The given $n$-dimensional faithful representation of $G$
is a subrepresentation of the regular representation,
since otherwise it would contain some irreducible subrepresentation with
multiplicity $>1$, which could be removed once to produce
a faithful subrepresentation on a lower-dimensional subspace,
contradicting the fact that the function $\dmax(n)$ is strictly increasing.
(Alternatively, this could be deduced from the fact that 
the given representation is irreducible for all $n\ne 9,10$, 
and is a direct sum of distinct irreducible
representations for $n=9$ and $n=10$.)
Moreover, Corollary~\ref{C:inverseGalois} shows that 
$G$ is realizable  as a Galois group over $\Q$.
Thus Proposition~\ref{P:subregular} yields $\al \in \Qbar$
with $n(\al)=n$ and $d(\al)=\#G=\dmax(n)$.
\end{proof}

\subsection{Explicit numbers attaining $\dmax(n)$}
\label{exceptionaln}

In theory, given $n \geq 1$, 
we can construct explicit $\al \in \Qbar$
with $n(\al)=n$ and $d(\al)=\dmax(n)$ as follows.
Let $G$\/ be a maximal-order finite subgroup of $\GL_n(\Q)$.
Take $e_j$ to be the column vector in $\Z^n$
having $j$-th entry $1$ and the rest $0$,
let $G_1$ be the stabilizer of $e_1$ under the left action of~$G$,
and put $N=|G:G_1|$, the size of the orbit of $e_1$ under this action.
For most of the groups we consider,
all of $e_1\cc e_n$ are in this orbit,
and so we denote the whole orbit by $e_1\cc e_n\cc e_N$.
We then find an {\it auxiliary polynomial}\/ $P_N$ of degree~$N$,
irreducible over $\Q$, whose splitting field has Galois group
$G$\/ over~$\Q$. Further, $n$ zeros $\be_1\cc \be_n$ of $P_N$
can be chosen so that the full list of conjugates
$\be_1\cc \be_N $ of $\be_1$ are the $(\be_1\cc \be_n)e_j$ for $j=1\cc N$.
 
The auxiliary polynomial $P_N$ arises, at least generically, as follows:
by Proposition~\ref{P:invariants},
we can write $\Q(x_1,\dots,x_n)^G=\Q(I_1,\dots,I_n)$,
where the $I_j$ are $G$-invariant homogeneous polynomials in the~$x_i$.
Choose $c_1,\dots,c_n \in \Q$,
and define a zero-dimensional variety $\mathcal V$\/ 
by the polynomial equations 
\begin{align*}
	I_1(x_1,\dots,x_n) &= c_1, \\
	&\vdots \\
	I_n(x_1,\dots,x_n) &= c_n. \\
\end{align*}
Then 
successively eliminate $x_n$, $x_{n-1}$, \dots, $x_2$
to get a monic polynomial $R(x_1)$ of degree $d_R$ given by $d_R=\prod_{j=1}^n\deg I_j$. 
Clearly $\mathbf xg \in \mathcal V$
for any $\mathbf x \in \mathcal V$\/ and $g\in G$, 
so the multiset of zeros of $R$ is 
$\{\mathbf xge_1\mid g\in G\}$, 
which consists of $\#G_1$ copies of
$\{\mathbf xe_j\mid j=1\cc N\}$. 
Thus $R(x_1)=P_N(x_1)^{\# G_1}$ for some polynomial $P_N$. 
For reflection groups and unitary reflection groups
we can choose the $I_j$ so that $d_R=\# G$;
in this case $P_N$ has degree $N$.
The polynomial $P_N$ is our auxiliary polynomial.

Choose $b_1,\dots,b_n \in \Q$ such that $b_1 x_1 + \cdots + b_n x_n$
is not fixed by any $g \in G$\/ except the identity. Then 
$\al=b_1\be_1+\cdots + b_n\be_n$ has $n(\al)=n$ and degree $\dmax(n)$,
its conjugates being $(\be_1 \cc \be_n)g(b_1 \cc b_n)^T$ for $g\in G$.
(This is the standard ``primitive element'' construction
for the Galois closure of $\Q(\be)$.)
For most choices of $(c_1,\dots,c_n)$
(that is, for all choices outside a ``thin set'', 
in the sense of~\cite{serretopics}), 
this construction will produce the required $\al$.
For small $n$ 
(such as $n=2$, considered in Sections \ref{n2ex} and~\ref{cyc_res}),
this procedure works well. 
For much larger $n$, however, the elimination process becomes impractical. 
Also, it becomes hard to check whether a particular
choice of $(c_1,\dots,c_n)$ yields a suitable~$\al$.
The difficulty is to choose $c_1\cc c_n$
so that not only is $P_N$ irreducible,
but also it has Galois group $G$\/
(instead of a subgroup).
For this reason, the following sections discuss more practical ways
of constructing $\al$, in the nonexceptional case and for $n=4$.

For the larger exceptional values of $n$,
even these methods would require special treatment for each value, 
and the large size of $\#G$\/ (see Table~\ref{Ta1}) 
has dissuaded us from trying to do the same for these~$n$.
One approach to constructing $\al\in\Qbar$ attaining $\dmax(n)$
for $6 \leq n \leq 10$ is to start
with Shioda's beautiful analysis relating the Weyl groups
of $E_6$,~$E_7$,~$E_8$ and their invariant rings
with the \MW\ lattices of rational elliptic surfaces
with an additive fiber.
For instance, in \cite[p.484--5]{shioda1991} Shioda uses this theory
to exhibit a monic polynomial in $\Z[X]$ with Galois group $W(E_7)$,
whose roots are the images of the 56 minimal vectors of the $E_7^*$
lattice under a \hbox{$\Q$\/-linear}, \hbox{$W(E_7)$-equivariant} map
{}from $E_7^*\otimes\Q$ to $\Qbar$.  The image under this map
of any vector in $E_7^*\otimes\Q$ with trivial stabilizer in $W(E_7)$
(that is, in the interior of a Weyl chamber) is then an $\al\in\Qbar$
with $n(\al)=7$ and $d(\al)=\#W(E_7)=\dmax(7)$.
A similar construction will work for $n=8$,
and (combined with the analysis of algebraic numbers
of conjugate dimension $1,2$) also for $n=9,10$.
The case $n=6$ will require additional work,
because Shioda's construction, which yields Galois group $W(E_6)$,
will have to be modified to produce $\langle W(E_6),-I\rangle$.

\subsection{Explicit numbers attaining $\dmax(n)$ for nonexceptional $n$}
\label{S:nonexceptional}

\begin{proposition}
\label{P:squarerootfield}
Let $k$ be a field of characteristic not~$2$.
Let $n \geq 2$.
Suppose $f(x)=x^n - a_1 x^{n-1} + \dots + (-1)^n a_n \in k[x]$ 
is a separable polynomial of degree $n$ with Galois group $S_n$
and discriminant $\Delta$.
Let $r_1,\dots,r_n \in \kbar$ be the zeros of $f(x)$.
Choose a square root $\sqrt{r_i}$ of each $r_i$, 
and let $K=k(\sqrt{r_1},\dots,\sqrt{r_n})$.
If $a_n \notin \Delta^\Z k^{*2}$ and 
either $n$ is even or $r_1 \notin k^* k(r_1)^{*2}$,
then $[K:k]=2^n n!$.
\end{proposition}

\begin{proof}
The action of the group $G:=\Gal(K/k)$ on
$\{\sqrt{r_1},-\sqrt{r_1},\dots,\sqrt{r_n},-\sqrt{r_n}\}$
is faithful and preserves the partition
$\{\{\sqrt{r_1},-\sqrt{r_1}\},\dots,\{\sqrt{r_n},-\sqrt{r_n}\}\}$,
so $G$\/ is a subgroup of the signed permutation group $W(B_n)$.
Recall that $W(B_n)$ is a semidirect product
	$$0 \to V \to W(B_n) \to S_n \to 1$$
where $V$\/ as a group with $S_n$-action 
is the standard permutation representation of $S_n$ over~$\F_2$.
Since $f$ has Galois group $S_n$,
the group~$G$\/ surjects onto the quotient $S_n$ of $W(B_n)$.
Considering the conjugation action of $G$\/ on itself 
gives a (possibly nonsplit) exact sequence
	$$0 \to W \to G \to S_n \to 1$$
for some subrepresentation $W$ of $V$.
The only subrepresentations of $V$\/ are $0$,
$\F_2$ with trivial $S_n$-action,
the sum-zero subspace of $V=\F_2^n$,
and $V$\/ itself.
If $W=V$, we are done.

If $W$\/ is contained in the sum-zero subspace,
then $W$\/ acts trivially on 
the square root 
$\be:=\sqrt{r_1} \dots \sqrt{r_n}$
of $a_n$.
Hence the action of $G$\/ on $\be$ is given by either the trivial
character or the sign character of $S_n$.
Thus either $\be \in k$ or $\be \sqrt{\Delta} \in k$.
Squaring yields $a_n \in \Delta^\Z k^{*2}$,
contrary to assumption.

The only remaining case is where $n$ is odd and $W=\F_2$.
Then $W$\/ acts trivially on 
the square root 
$\be_1:=\sqrt{r_2} \sqrt{r_3} \dots \sqrt{r_n}$
of $r_2 r_3 \dots r_n = a_n/r_1$.
Hence the action of $\Gal(K/k(r_1))$ on $\be_1$ is given by 
either the trivial character or the sign character 
of $S_{n-1} = \Gal(k(r_1,\dots,r_n)/k(r_1))$.
Thus either $\be_1 \in k(r_1)$ or $\be_1 \sqrt{\Delta} \in k(r_1)$.
Squaring shows that $r_1 \in k^* k(r_1)^{*2}$,
again contrary to assumption.
\end{proof}

In the situation of Proposition~\ref{P:squarerootfield},
when its hypotheses are satisfied,
we can take the auxiliary polynomial to be $P_{2n}(x)=f(x^2)$.

The following corollary is needed in Section~\ref{n=4}.

\begin{corollary}
\label{C:totallyreal}
Let  $n \geq 2$.
Suppose $f(x)=x^n - a_1 x^{n-1} + \dots + (-1)^n a_n \in k[x]$ 
is a polynomial of degree $n$ over a field $k\subset \R$, with Galois group $S_n$.
Suppose that the zeros $r_1,\dots,r_n$ of $f(x)$ are real and
satisfy $r_1 < 0 < r_2 < \dots < r_n$.
Choose a square root $\sqrt{r_i} \in \overline k$ of each $r_i$.
and let $K=k(\sqrt{r_1},\dots,\sqrt{r_n})$.
Then $[K:k]=2^n n!$.
\end{corollary}

\begin{proof}
It suffices to check the hypotheses
of Proposition~\ref{P:squarerootfield}.
The discriminant $\Delta$ satisfies $\Delta>0$,
but $a_n = r_1 \dots r_n < 0$, so $a_n \notin \Delta^\Z k^{*2}$.

If $r_1 \in k^* k(r_1)^{*2}$,
say $r_1 = c \ga_1^2$ with $c \in k^*$ and $\ga_1 \in k(r_1)$,
then applying an automorphism yields $r_2 = c \ga_2^2$
with $\ga_2 \in k(r_2)$.
These two equations force $c<0$ and $c>0$, respectively,
a contradiction.
\end{proof}

\begin{proposition}
\label{P:nonexceptionalexample}
For $n=1$ let $r_1=2$, while for $n \geq 2$
let $r_1,\dots,r_n \in \Qbar$ be the zeros of $f(x)=x^n + (-1)^n(x-1)$.
Choose a square root of each $r_i$, 
and let $\al=\sqrt{r_1} + 2\sqrt{r_2} + \dots + n\sqrt{r_n}$.
Then $n(\al)=n$ and $d(\al)=2^n n!$.
\end{proposition}

\begin{proof}
By~\cite[p.~42]{serretopics}, 
the polynomial $(-1)^n f(-x) = x^n-x-1$ has Galois group $S_n$ over $\Q$,
so $f(x)$ has Galois group $S_n$ over $\Q$.
Also by~\cite[p.~42]{serretopics},
each inertia group of $\Gal(\Q(r_1,\dots,r_n)/\Q)$
is either trivial or generated by a transposition;
it follows that the same is true for the Galois group
$G$\/ of $f$\/ over $\Q(i)$.
The group $G$\/ has index at most $2$ in $S_n$,
so $G$\/ is $S_n$ or $A_n$.
We claim that $G=S_n$.
For $n=2$ we check this directly.

Take $n \geq 3$. If $G=A_n$, then
as $G$\/ would contain no transpositions,
 all the inertia groups in $G$\/ would be trivial,
and $\Q(i)$ would have an $A_n$-extension unramified
at all places.
The existence of such an extension contradicts 
the Minkowski discriminant bound for $n \geq 4$,
and contradicts class field theory for $3 \leq n \leq 4$.
Thus $G=S_n$.

In particular, if $\Delta$ is the discriminant of $f(x)$,
then $\Delta \notin \Q(i)^{*2}$, so $|\Delta| \notin \Q^{*2}$.
Therefore $a_n := -1$ is not in $\Delta^\Z \Q^{*2}$.

We now finish checking the hypotheses in Proposition~\ref{P:squarerootfield}
by showing that the assumptions $n$ odd and $r_1 \in \Q^* \Q(r_1)^{*2}$
lead to a contradiction.
Suppose $n$ is odd,
and $r_1 = c \ga^2$, with $c \in \Q^*$ and $\ga \in \Q(r_1)^*$.
Taking $N_{\Q(r_1)/\Q}$ of both sides yields 
$(-1)^n \equiv c^n \pmod{\Q^{*2}}$.
Since $n$ is odd, $c \equiv -1 \pmod{\Q^{*2}}$.
Without loss of generality, $c=-1$.
Since $\ga$ generates $\Q(r_1)$, 
the monic minimal polynomial $g(t) \in \Q[t]$ of $\ga$
is of degree $n$.
Write $g(t)g(-t) = h(t^2)$ for some polynomial $h \in \Q[x]$.
Substituting $t=\ga$ shows that $h(-r_1)=0$, 
but $h$ has degree $n$, so $h(x)=f(-x)$.
Thus the polynomial $-f(-t^2) = t^{2n}-t^2-1$ factors as $- g(t) g(-t)$.
However, it is known to be irreducible (Ljunggren \cite[Theorem 3]{ljunggren1960}).

By Proposition~\ref{P:squarerootfield},
the field $K=\Q(\sqrt{r_1},\dots,\sqrt{r_n})$ has degree $2^n n!$.
Each $\sqrt{r_i}$ lies outside the field generated by the other
square roots over $\Q(r_1,\dots,r_n)$,
so $\sqrt{r_1}$, \dots, $\sqrt{r_n}$ are linearly independent over $\Q$.
The conjugates of $\al$ are the numbers of the form
$\sum_{j=1}^n \ep _j j \sqrt{r_{\sigma(j)}}$
where $\sigma \in S_n$ and $\ep _1,\dots,\ep _n \in \{\pm 1\}$.
The linear independence of the square roots guarantees 
that these $2^n n!$ elements are distinct.
\end{proof}

\subsection{An explicit number attaining $\dmax(n)$ for $n=2$}\label{n2ex}

For $n=2$, we can take $P_6(x)=x^6-2$. Taking one zero $\be$ of $P_6$,
all zeros are spanned by the two zeros $\be, \om_3 \be$
 where $\om_3$ is  a primitive
cube root of unity. Then
$\al=\be + 3 \om_3\be$
has $n(\al)=2$ and $d(\al)=12$, and minimal polynomial $y^{12}+572 y^6 + 470596$.

\begin{remark}
This example can be produced using the procedure outlined in Section~\ref{exceptionaln},
as follows.
The group $W(G_2)$ from Table~\ref{Ta1} equals 
$\langle \left( \begin{smallmatrix}
0 & -1 \\
1 & \phantom{-}1
\end{smallmatrix}
\right),
\left( \begin{smallmatrix}
0 & 1 \\
1 & 0
\end{smallmatrix}
\right)
\rangle$,
and has invariants $I_1=x_1^2-x_1x_2+x_2^2$ and $I_2=(x_1x_2(x_1-x_2))^2$.
Taking $c_1=0$, $c_2=2$, $b_1=1$, $b_2=-3$,
we get the minimal polynomial of $\al$ as the $x_2$-resultant of 
$I_1(y+3x_2,x_2)$ and $I_2(y+3x_2,x_2)-2$.
\end{remark}

\subsection{An explicit number attaining $\dmax(n)$ for $n=4$}
\label{n=4}

For $n=4$, one maximal-order finite subgroup of $\GL_4(\Q)$
is the order-$1152$ group $W(F_4)$ 
generated by its index-$3$ subgroup $W(B_4)$ (of order $384$)
and the order-$2$ matrix
\[
\frac12	\begin{pmatrix} 
   1 & \phantom{-}1 & \phantom{-}1 & \phantom{-}1 \\
   1 & -1 & \phantom{-}1 & -1 \\
   1 & \phantom{-}1 & -1 & -1 \\
   1 & -1 & -1 & \phantom{-}1 \\
	   \end{pmatrix}.
\]
Thus by Galois correspondence we should be able to apply
the construction of Section~\ref{exceptionaln} 
for $\be$ defined over a suitable cubic extension of $\Q$. 
And indeed, this is possible. 
  
Define $s_{2k}=z_1^{2k}+z_2^{2k}+z_3^{2k}+z_4^{2k}$ for $k=1,2,\dots $.
Four independent homogeneous invariants for $W(F_4)$ are known \cite{mehta1988} to be
$$ 
I_{2k}= (8-2^{2k-1})s_{2k}+\sum_{j=1}^{k-1}\binom{2k}{2j}s_{2j}s_{2k-2j}
$$
for $k=1,3,4,6$. Using the Newton identities and with the help of Maple
these can be written entirely as polynomials
in $s_2,s_4,s_6,s_8$ as follows:
\begin{gather*}
 I_2=6s_2, \qquad 
{I_{6}} =  - 24\,{s_{6}} + 30\,{s_{2}}\,{s_{4}}, \qquad
{I_{8}}=  - 120\,{s_{8}} + 56\,{s_{2}}\,{s_{6}} + 70\,s_4^2, \\\label{E:Invs}
 \begin{split}
{I_{12}}  =  - 540\,{s_{4}}\,{s_{8}} + 244\,s_6^2 - 1365
\,s_2^2\,{s_{8}} + {\displaystyle \frac {1365}{2}} \,s_2
^2\,s_4^2 + 255\,s_4^3\\
 - 710\,s_2^4\,{s_{4}} + 1250\,s_2^3 \,{s_{6}}
 + {\displaystyle \frac {159}{2}} \,s_2^6 + 110\,{s_{2}}\,{s_{4}}\,{s_{6}}. 
\end{split}
\end{gather*}

We now use resultants to eliminate $s_4$ and $s_6$. This shows that $s_8$ is
cubic over $\Q(I_2,I_6,I_8,I_{12})$, and also that $s_4,s_6\in \Q(I_2,I_6,I_8,I_{12})(s_8)$.
Specifically, we take $I_2=6s_2 = 30,I_6=1410,I_8=13670$
and $I_{12}=1161749$, and then
$\ga:=s_8$ (the real root, say) satisfies
\[
 \ga^{3} + {\displaystyle \frac {5735}{32}} \,\ga^{2}
 + {\displaystyle \frac {5811288377}{36864}} \,\ga - 
{\displaystyle \frac {114051068048293}{6220800}} 
=0. 
\]
  Then, with the Newton identities,
we compute the values of the elementary symmetric functions
of the $z_i^2$. This gives a polynomial $Q_4$ satisfied by the $z_i^2$:
{\small
\begin{multline*}\label{E:minpoly4_be}
Q_4(x)=\ x^{4} - 5\,x^{3} + {\displaystyle \frac {
20261200695}{3175710433}} \,x^{2} + {\displaystyle \frac {34560}{
3175710433}} \,x^{2}\,\ga^{2} - {\displaystyle \frac {
47690820}{3175710433}} \,x^{2}\,\ga \\
 + {\displaystyle \frac {36679035170}{9527131299}} \,x - 
{\displaystyle \frac {28800}{3175710433}} \,x\,\ga^{2} + 
{\displaystyle \frac {39742350}{3175710433}} \,x\,\ga - 
{\displaystyle \frac {203476507483}{38108525196}}  \\
 - {\displaystyle \frac {72000}{3175710433}} \,\ga^{2}
 - {\displaystyle \frac {56249419}{12702841732}} \,\ga .
\end{multline*}
}

 We write its zeros as $\be^2_1,\be^2_2,\be^2_3,\be^2_4$ say.
 They are real and close to $-1,1,2$, and $3$.
(The values for the invariants were chosen to be close
to the values they would have had if $z_i^2, i=1 \cc 4$  
had been {\it exactly} $-1,1,2,3$.) Furthermore,  its discriminant 
 $223967999/97200 $ 
 is not a square in $\Q(\ga)$.  Now, shifting $x$
in  this quartic by $5/4$ to obtain a polynomial 
$z^4+b_2z^2+b_1z+b_0$ having zero cubic term,
its cubic resolvent $z^3+2b_2z^2+(b_2^2-4b_0)z-b_1^2$
is readily checked to be irreducible over $\Q(\ga)$.
Hence by \cite[Ex.~14.7, p. 117]{garling1986},
the Galois closure of $\Q(\ga,\be)$ over $\Q(\ga)$
has Galois group $S_4$. Then, as $\be^2_1<0< \be^2_2<\be^2_3<\be^2_4$,
we have
$[\Q({\be_1},{\be_2},{\be_3},{\be_4}):\Q]
= 2^4\cdot 4! = 384$, 
on applying Corollary~\ref{C:totallyreal} with $k=\Q(\ga)$. 

If we now take the resultant of $Q_4(x^2)$ 
and the minimal polynomial of~$\ga$, to eliminate $\ga$, 
we obtain the degree $24$ auxiliary polynomial
{\small
\begin{multline*}
P_{24}(x)= \\
x^{24} - 15\,x^{22} + {\displaystyle \frac {
375}{4}} \,x^{20} - {\displaystyle \frac {2405}{8}} \,x^{18} + 
{\displaystyle \frac {65435}{128}} \,x^{16} - {\displaystyle 
\frac {25905}{64}} \,x^{14} - {\displaystyle \frac {181583}{3072}
} \,x^{12} + {\displaystyle \frac {8367137}{18432}} \,x^{10} \\
 - {\displaystyle \frac {28198575}{65536}} \,x^{8} + 
{\displaystyle \frac {1338226651}{5308416}} \,x^{6} - 
{\displaystyle \frac {895964239}{8847360}} \,x^{4} + 
{\displaystyle \frac {4234139}{294912}} \,x^2 - {\displaystyle 
\frac {24389830879}{1592524800}}.
\end{multline*}
}

 This polynomial is irreducible, with zeros
$\frac{1}{2}( \pm{\be_1}\pm {\be_2}\pm {\be_3}\pm {\be_4})$
as well as $\pm{\be_1}$, $\pm {\be_2}$, $\pm {\be_3}$, $\pm {\be_4}$.
Now $(1,2,3,5)^T$ is not a fixed point of any $g\ne I$ in $W(F_4)$. It follows that $
\al = {\be_1} +2 {\be_2} +3 {\be_3} +5 {\be_4}$
has $n(\al)=4$ and degree $d(\al)= 1152$, its conjugates being the numbers 
$({\be_1}, {\be_2},{\be_3}, {\be_4})g(1,2,3,5)^T$ for
$g\in W(F_4)$.

\section{Conjugate dimensions over other fields}
\subsection{General results}
\label{gen_res}
The conjugate dimension
can behave differently if we use ground fields other than~$\Q$.
For a field~$k$ and a positive integer~$n$,
let $D(k,n)$ be the maximal degree of $\al\in k^s$
of $k$-conjugate dimension at most~$n$.  For instance $D(\Q,n) = \dmax(n)$.
If the degree is unbounded, we set $D(k,n)=\infty$.
This can happen even for Hilbertian fields of characteristic zero. 
For example, $D(\C(t),1) = \infty$, because for each $d \geq 1$
a \hbox{$d$\/-th} root of~$t$ generates the Galois extension
$\C(t^{1/d})$ of degree~$d$, and all conjugates of $t^{1/d}$
generate the same \hbox{$1$-dimensional} space.
Nevertheless we can generalize some of our results
to various ground fields other than~$\Q$.
We obtain the following.

\begin{theorem}\label{T:otherfields}
$\left.\right.$
\begin{enumerate}
\item[(i)]
If $k$ is a number field of degree $m$ over $\Q$,
then $\dmax(n)\leq D(k,n) \leq \dmax(mn)$ for all $n \geq 1$.
\item[(ii)]
If $k$ is a Hilbertian field of characteristic not dividing $\l$
and $k$ contains the $ \l$-th
 roots of unity, then $D(k,n) \geq \l^n n!$.
\item[(iii)]
If $k$ is a finitely generated transcendental extension of $\C$,
then $D(k,n)=\infty$ for all $n \geq 1$.
\item[(iv)]
If $k$ is a finite field of $q$ elements, then $D(k,n) = q^n-1$.
\item[(v)]
If $k$ is a finitely generated transcendental extension 
of a finite field $k_0$, then $D(k,1)=q-1$ where $q$ is the size 
of the largest finite subfield of $k$,
and $D(k,n)=\infty$ for all $n \geq 2$.
\end{enumerate}
\end{theorem}

\begin{proof}
$\left.\right.$

(i)
By Proposition~\ref{P:faithful}, if $\al \in k^s$
has degree~$d$ and conjugate dimension~$n$
then there exists a \hbox{$d$\/-element} subgroup of $\GL_n(k)$.
If $[k:\Q]=m$, then an $n$-dimensional vector space over $k$
can be viewed as an $mn$-dimensional vector space over $\Q$,
so we get an injection $\GL_n(k) \hookrightarrow \GL_{mn}(\Q)$.
Hence $d \leq \dmax(mn)$.  
For the lower bound, note that the specialization made in
Proposition~\ref{P:hilbert} can, by \cite[Theorem~46, p.~298]{schinzel2000}, 
be made in such a way that
the minimal polynomial of the algebraic number 
with conjugate dimension $n$ remains irreducible over
the field $k$. 
This gives an example of an algebraic number of degree $\dmax(n)$ 
over $k$ and $k$-conjugate dimension at most $n$, 
so $\dmax(n)\leq D(k,n) $.

(ii)
If $k$ contains the $\l$-th roots of unity then
$\GL_n(k)$ contains the group of size $\l^n n!$
consisting of the permutation matrices 
whose entries are $\l$-th roots of unity.
Moreover, the invariant ring of this group is polynomial,
being generated by the elementary symmetric functions
of the \hbox{$\l$-th} powers of the coordinates.
Thus the invariant field is purely transcendental over~$k$.
Therefore, by Propositions \ref{P:subregular} and~\ref{P:hilbert},
there exist $\al \in k^s$ of conjugate dimension~$n$ and degree $\l^n n!$.

(iii) 
This follows from (ii), using  the fact that every such
field is Hilbertian (\cite[ Theorem 49, p. 308]{schinzel2000}).

(iv)
The Galois group of any $k(\al)/k$ with $n(\al)=n$ must be contained in $\GL_n(k)$,
but must also be cyclic because $k$ is a finite field $\F_q$.
Hence $\# G \leq q^n-1$, as may be seen using the characteristic equation
of an invertible matrix in $\GL_n(k)$.  We claim that
the field of $q^{q^n-1}$ elements is generated by an element~$\al$
of conjugate dimension~$n$ over~$k$.
Let $g$ be a generator of $\F_{q^n}^*$,
and let $f(x)=\sum_{i=0}^n c_i x^i$ be its minimal polynomial over $\F_q$.
Let $\al \in \Fbar_q^*$ be a zero of $\sum_{i=0}^n c_i X^{q^i}$.
Make the $\F_q$-vector space $\Fbar_q$ 
into a module over the polynomial ring $\F_q[\tau]$
by letting $\tau$ act as the endomorphism $z \mapsto z^q$.
Then the ideal $I$ of $\F_q[\tau]$ 
that annihilates $\al$ contains $f(\tau)$,
but $I \ne (1)$.
Since $f$\/ is irreducible, $I=(f(\tau))$.
Thus the \hbox{$\F_q$-span} of $\al$ and its conjugates 
is an $\F_q[\tau]$-module isomorphic to $\F_q[\tau]/(f(\tau))$.
In particular, $n(\al)=\deg f =n$.
Also $d(\al)$ is the smallest $d$ such that $\tau^d(\al)=\al$,
which is the smallest $d$ such that $\tau^d=1$ in $\F_q[\tau]/(f(\tau))$;
by choice of $g$, we get $d=q^n-1$.

(v)
Without loss of generality, suppose that $k_0$
is the largest finite subfield 
  of $k$, so $\# k_0=q$.
Suppose $\al \in \kbar$ has $n(\al)=1$.
Proposition~\ref{P:faithful} bounds $d(\al)$
by the size of the largest finite subgroup of $\GL_1(k)=k^*$.
Elements of finite order in $k^*$ are roots of unity,
hence contained in $k_0^*$.
Thus $D(k,1) \leq q-1$.
The opposite inequality follows from (ii) since,
by \cite[Theorem 47, p. 301]{schinzel2000}, $k$ is Hilbertian.

Now suppose $n \geq 2$.
Choose a finite Galois extension $L$ of $k$ with $[L:k]=n-1$.
(For instance, let $L$ be the compositum of a suitable subfield of a
cyclotomic extension of $k$ 
with some Artin-Schreier extensions of $k$.)
Let $V$\/ be the $\F_q$-span
of a $\Gal(L/k)$-stable finite subset of~$L$
that spans $L$ as a $k$-vector space.
Define
$$
P_{V,\ep }(X) := \prod_{x\in V} (X-x) + \ep 
\in k[X,\ep ],
$$
where $\ep$ is an indeterminate.
Then $P_{V,0}(X)$ is a $q$-linearized polynomial in~$X$, that is,
a \hbox{$k$-linear} combination of $X,X^q,X^{q^2},\ldots$.
(See \cite[Corollary~1.2.2]{goss1996}, for instance.)
It has distinct roots, namely the elements of~$V$.
Therefore $P_{V,\ep }(X)$, considered as a polynomial in~$X$,
has distinct roots, which constitute a translate of~$V$\/
in the separable closure of~$k(\ep )$.
Moreover, $P_{V,\ep }(X)$ is irreducible,
because it is a monic polynomial in~$\ep $ of degree~$1$.
Since $k$\/ is Hilbertian, it contains $c \ne 0$ such that
$P_{V,c} \in k[X]$ is irreducible.
Let $\al$ be a zero of~$P_{V,c}$.
Then $\al$ is an element of $k^s$ of degree~$\# V$.
Since the set of conjugates of $\al$ is $\{\al+v\mid v\in V\}$, 
the $k$-span of this set
equals the span of $V \cup \{\al\}$.
However $\al \notin L$ 
since $d(\al) = \#V \geq q^{n-1} > n-1$.
So, as the $k$-span of $V$\/ is $L$,  $n(\al)=[L:k]+1=n$.
Thus $D(k,n) \geq \#V$.
Since $V$\/ can be taken arbitrarily large, $D(k,n)=\infty$.
\end{proof}

\subsection{Results for cyclotomic fields}
\label{cyc_res}

Theorem~\ref{T:main} generalizes to finite cyclotomic extensions of~$\Q$.
Let $\om_\l$ be a primitive $\l$-th root of unity.

\begin{theorem}\label{T:cyc}
Fix an integer $n \geq 0$ and an even integer $\l\geq 4$.
If $\al \in \Qbar$ has conjugate dimension $n$ over $\Ql$
then the degree $d$ of $\al$ over $\Ql$ satisfies
\[
	n \leq  d \leq D(\Ql,n),
\]
where $D(\Ql,n)$ is defined by Table~\ref{Ta2}.
In particular, $D(\Ql,n)=\l^n n!$ for 
\[
	(n,\l)\notin 
	\{(2,4),(2,8),(2,10),(2,20),(4,4),(4,6),(4,10),(5,4),(6,6),(6,10),(8,4)\}.
\]
Furthermore, for each pair $(n,\l)$ with $n \geq 1$ and $\l\geq 4$ even,
there exist $\al \in \Qbar$ attaining the lower and upper bounds.
\end{theorem}

\begin{table}[ht]
\begin{center}
\bigskip\begin{tabular} {r | r | r  | l | r}

 $n$  & $\ell$ & $D(\Q(w_\ell),n) / (\ell^n n!) $ & Maximal-order subgroup $G$
  & $D(\Ql,n)=\# G$\\
     \hline
    2 & 4 & 3   & $\ST_8 =\langle \GL_2(\F_3),\om_4I\rangle$& 96\\
    2 & 8 & 3/2 & $\ST_9=\langle \GL_2(\F_3),\om_8I\rangle $& 192\\
    2 &10 & 3   & $ \ST_{16}=\langle \om_5I\rangle\times \SL_2(\F_5) $ & 600\\ 
    2 &20 & 3/2 & $\ST_{17}=\langle \SL_2(\F_5),\om_{20}I\rangle$ & 1200\\ 
    4 & 4 & 15/2& $\ST_{31}$                          & 46080\\
    4 & 6 & 5   & $\ST_{32} $                         & 155520\\
    4 &10 & 3   &  $ \ST_{16}\wr S_2$ & 720000\\
    5 & 4 & 3/2 & $\ST_{31}\times \langle \om_4 I\rangle $ & 184320\\
    6 & 6 & 7/6 &  $\ST_{34}$  & 39191040\\
    6 &10 & 9/5 &  $ \ST_{16}\wr S_3$ & 1296000000\\
    8 & 4 &45/28&  $\ST_{31}\wr S_2$ & 4246732800\\
\multicolumn{3}{l|}{all  other $(n,\ell)$, $\l\geq 4$  even \quad\hfill 1}
& $ \ST_2(\ell,1,n) = (\Z/\ell\Z)^n \rtimes S_n \! $
& $\ell^n n!$ \\
    \end{tabular}    
    \vspace*{1ex}
    \caption{Maximal-order subgroups of $\GL_n(\Q(\om_\ell))$ for $\l\geq 4$ even }\label{Ta2}
\end{center}\end{table}

Table~\ref{Ta2} is a list of groups isomorphic to 
maximal-order finite subgroups $G$\/ of $\GL_n(\Ql)$, 
quoted from Feit \cite{feit1996}. 
(An error in the first line of his table has been corrected.) 
In this table $\ST_j$ refers to the $j$-th 
unitary reflection group in \cite[Table VII]{shephard-todd1954}, 
and the wreath product $G \wr S_n$ is the semidirect product
$(G \times \cdots \times G) \rtimes S_n$ 
in which $S_n$ acts on the $n$-fold product of $G$\/ 
by permuting the coordinates.
See also \cite[Table 7.3.1]{smith1995}.

\begin{proof}
The proof is a generalization of that of Theorem~\ref{T:main}.
For fixed $\l$, $D(\Ql,n)$ is a strictly increasing function of~$n$.
Thus to carry over the proof, it remains to show that the invariant
subfield $\Ql(x_1,\dots,x_n)^G$ is purely transcendental over $\Ql$
in each case of Table~\ref{Ta2}.
This is immediate for all the Shephard-Todd groups in the table,
by the extension of Chevalley's Theorem
to unitary reflection groups by Shephard and Todd (\cite{shephard-todd1954};
see also \cite[p.~115,~Thm.~4]{bourbaki-lie456}, \cite[p.~65]{humphreys1990}).
For example, when $G=(\Z/\ell\Z)^n \rtimes S_n$,
the field of invariants $\Ql(x_1\cc x_n)^G$ is $\Ql(e_1\cc e_n)$, 
where $e_j$ is the $j$-th elementary symmetric
function of $x_1^\ell\cc x_n^\ell$.
The three remaining cases are handled
by Lemma~\ref{L:wreath invariants} below.
\end{proof}

\begin{lemma}
\label{L:diagonal action}
Let $k$ be a field.
Let the symmetric group $S_m$ act on
\[
	K=k(x_1^{(1)},\dots,x_1^{(m)}; \dots; x_n^{(1)},\dots,x_n^{(m)})
\]
by acting on the superscripts.
Then $K^{S_m}$ is purely transcendental over~$k$.
\end{lemma}

\begin{proof}
If $E/F$ is a Galois extension of fields with Galois group $G$, and
$V$\/ is an $E$-vector space equipped with a semilinear action of~$G$,
there exists an $E$-basis of $V$\/ consisting of $G$-invariant 
vectors~\cite[II.5.8.1]{silvermanAEC}.

Apply this to $E=k(x_1^{(1)},\dots,x_1^{(m)})$, $G=S_m$, 
$F=E^G$ (the purely transcendental extension of $k$
generated by the symmetric functions in $x_1^{(1)},\dots,x_1^{(m)}$),
and $V$\/ the \hbox{$E$-subspace} of~$K$\/
spanned by all the $x_i^{(j)}$ with $i \ge 2$.
Choose an $E$-basis $\{v_s\}$ of $G$-invariant vectors as above.
Let $K_0=k(\{v_s\})$.
Since $EK_0=K$, we have $[K:K_0] \le [E:F]=m!$,
On the other hand, $K_0 \subseteq K^G$ with $[K:K^G]=m!$,
so $K_0=K^G$.
Since the $x_i^{(j)}$ are algebraically independent over $E$,
the $v_s$ are algebraically independent over $k$.
\end{proof}

\begin{lemma}  
\label{L:wreath invariants}
Let $k$ be a field,
and let $G$ be a finite subgroup of $\GL_n(k)$ 
whose field of invariants $k(x_1,\dots,x_n)^G$
is purely transcendental over $k$.
Let $G \wr S_m$ act on
\[
	L=k(x_1^{(1)},\dots,x_n^{(1)}; \dots; x_1^{(m)},\dots,x_n^{(m)})
\]
by letting the $i$-th of the $m$ copies of $G$ act linearly
on the span of $x_1^{(i)},\dots,x_n^{(i)}$
while $S_m$ acts on the superscripts.
Then $L^{G \wr S_m}$ is purely transcendental over~$k$.
\end{lemma}

\begin{proof}
Since $G \wr S_m$ is a semidirect product of $S_m$ by $G^m$,
we have $L^{G \wr S_m} = \left(L^{G^m} \right)^{S_m}$.
If $k(x_1,\dots,x_n)^G=k(I_1,\dots,I_n)$,
then
\[
	L^{G^m} = k(I_1^{(1)},\dots,I_n^{(1)}; \dots;
			I_1^{(m)},\dots,I_n^{(m)}),
\]
and $S_m$ acts on this by acting on superscripts.
Now apply Lemma~\ref{L:diagonal action}.
\end{proof}

{\it Example.} Using the elimination procedure outlined in
Section~\ref{exceptionaln},
we can give an example of an algebraic number $\al$ of degree $96$ over $\Q(i)$
with $\Q(i)$-conjugate dimension $2$ and Galois group 
$\ST_8$, as in Table~\ref{Ta2}. Now 
$\ST_8= \langle \left( \begin{smallmatrix}
0 & 1 \\
1 & i 
\end{smallmatrix} \right), 
\left( \begin{smallmatrix}
0 & 1 \\
 - i & 0
\end{smallmatrix} \right) \rangle$, 
with invariants 
{\tiny
\begin{align*}
I_8(x_1,x_2) &= x_1^{8}+ 4(1 + \,i)\,x_1^{7}\,x_2+ 14\,i\,x_1^{6}\,x_2^{2} - 14(1 - \,i)\,
x_1^{5}\,x_2^{3}- 21\,x_1^{4}\,x_2^{4} - 14(1 + \,i)\,x_1^{3}\,x_2^{5}
 - 14\,i\,x_1^{2}\,x_2^{6}+ 4(1 - \,i)\,x_1\,x_2^{7} 
+ x_2^{8}      \\
I_{12}(x_1,x_2) &=  2\,x_1^{12}+ 12(1 + \,i)\,x_1^{11}\,x_2 
+ 66\,i\,x_1^{10}\,x_2^{2} 
- 110(1 - \,i)\,x_1^{9}\,x_2^{3}- 231\,x_1^{8}\,x_2^{4} 
\\
 &\qquad {} - 132(1 + \,i)\,x_1^{7}\,x_2
^{5} 
- 132(1 - \,i)\,x_1^{5}\,x_2^{7}
- 231\,x_1^{4}\,x_2^{8} 
- 110(1 + \,i)\,x_1^{3}\,x_2^{9} 
- 66\,i\,x_1^{2}\,x_2^{10}
+ 12(1 - \,i)\,x_1\,x_2^{11}
+2\,x_2^{12}.
\end{align*}
}
The $x_2$-resultant of $I_8-1-i$  and $I_{12}-1$ is 
$P_{24}(x_1)^4$, where the auxiliary polynomial $P_{24}$ is 
\[
P_{24}(x)= 27\,x^{24} - 270(1+i)\,x^{16}  + 270\,
x^{12} - 810\,i\,x^{8} + 54(1+i)\,x^{4}  - 9 + 8\,i.
\]
Two zeros $\be$ and $\be'$ of $P_{24}$ can be chosen so that the conjugates of
$\be$ are 
$$
\om \be,\quad \om \be',\quad \om (\be~+~\be'),\quad \om (\be-i\be'), \quad\om (\be+(1-i)\be'), 
\quad \om ((1+i)\be+\be'),
$$
 where $\om\in\{\pm 1,\pm i\}$. Then $\al=\be+2\be'$ has degree $96$ over
$\Q(i)$, with
conjugates $(\be,\be')g (1, 2)^T$ for $g\in \ST_8$. 
The minimal polynomial 
of $\al$ can be computed directly as the $x_2$-resultant of $I_8(y-2x_2,x_2)-1-i$
and $I_{12}(y-2x_2,x_2)-1$.

\subsection{$D(k,n)$ depends on more than $\l$ and $n$}

Let $k$ be a number field, and let $\l$ be the number
of roots of unity in $k$.
It seems reasonable to guess, as in the case of cyclotomic fields $\Ql$,
that $D(k,n) = \l^n n!$ for all but finitely many $n$.
However, it is possible that two number fields $k$ and $k'$ 
contain the same number of roots of unity, 
but $D(k,n)\ne D(k',n)$ for some $n$. 
For example, we can take
$k=\Q(\cos(2\pi/m),\sin(2\pi/m))$, where $m>6$, and $k'=\Q$. 
In both cases $\l=2$, but $D(k,2)>D(\Q,2)=12$. 
Indeed, there exist $a,b\in k$ such that 
$\al=\sqrt[m]{a}(1+b\om_m)$ is of degree $2m>12$ over $k$. 
Its  conjugate dimension  over $k$ is $2$;
 its conjugates are
spanned by $\sqrt[m]{a}$ and $i\sqrt[m]{a}$.  
This example also shows that the number of exceptional cases 
can be arbitrarily large, since we may simply take $m$ with $2m>2^n n!$.

Another example is $D(\Q(\sqrt 5\,),3)\geq 120$,
obtained from the icosahedral subgroup of $\GL_3(\Q(\sqrt 5\,))$
(reflection group $\ST_{23}$)
via Propositions \ref{P:subregular} and~\ref{P:hilbert}.

\section{Multiplicative conjugate rank}\label{multrank}

Instead of the dimension $n(\al)$ of the $\Q$-vector space spanned
by the $d$\/~conjugates $\al_i$ of an algebraic number~$\al$,
we may consider the rank~$r(\al)$
of the multiplicative subgroup of~$\Qbar^*$ they generate.
We call this the {\em (multiplicative) conjugate rank}\/ of~$\al$.
As before, we have the trivial inequality $r(\al) \leq d(\al)$,
which is sharp in the case of maximal Galois group
(again by~\cite[Lemma~1]{smyth1986}).  
Unlike in the additive case,
we can have no nontrivial lower bound without some further hypothesis,
because if $\al$ is a root of unity
then $r(\al)=0$ while $d(\al)$ is unbounded. However, also unlike the additive
case, we have the following result over a very general field. 
The main difficulty in the proof below is to show
that this bound is sharp for Hilbertian fields.

\begin{theorem}\label{T:mult}
Suppose that $\al$ is separable and 
algebraic of degree $d(\al)$ over a field~$k$, 
and the multiplicative subgroup of~$(k^s)^*$ generated
by the conjugates $\al_1,\ldots,\al_d$ of~$\al$ is torsion-free.
Then the rank $r(\al)$ of this subgroup satisfies
$r(\al) \leq d(\al) \leq \dmax(r(\al))$,
with $\dmax(\cdot)$ defined by Table~\ref{Ta1}
 as before.
If $k$\/ is Hilbertian, then for each integer $r \geq 1$
there are $\al\in k^s$ of conjugate rank~$r$
attaining the lower and upper bounds.
\end{theorem}

The upper bound is given by the same function $\dmax(\cdot)$
that we found for the conjugate dimension over $\Q$,
and this bound is independent of the ground field~$k$,
although it need not always be sharp.

\begin{proof}
For any $\al\in k^s$, let $\Gamma=\Gamma(\al)$
be the multiplicative group generated by the $\al_i$.
We observed already that the lower bound $d(\al) \geq r(\al)$ is immediate.
For the upper bound, we argue as we did for $n(\al)$.
The Galois group~$G$\/ acts faithfully on~$\Gamma$.
By hypothesis, $\Gamma \cong \Z^{r(\al)}$,
so $G$\/ acts faithfully also on $\Gamma \otimes_\Z \Q$,
which is a $\Q$-vector space of dimension $r(\al)$.
Hence $\#G$\/ is bounded above by $\dmax(r(\al))$,
the size of the largest finite subgroup of $\GL_{r(\al)}(\Q)$.
Hence $d(\al) \leq\# G \leq \dmax(r(\al))$.

The proof that there are examples attaining equality
when $k$ is Hilbertian uses two corollaries
of the following technical result.

\begin{proposition}
\label{P:multiplicative group contains ZG}
Let $L/k$ be a finite Galois extension of fields
with Galois group $G$, and suppose that $k$ is not algebraic over
a finite field.  
Then the ${\Z}G$-module $L^*$ contains a free ${\Z}G$-module of rank~$1$.
\end{proposition}

\begin{proof}
For each $g \in G - \{1\}$, choose $a_g \in L$ that is not fixed by $g$.
Choose $b \in L$ that is not algebraic over a finite field.
Let $S$ be the union of the $G$-orbits of the $a_g$ and of $b$.
Then $S$ is finite.
Let $L_0$ be the minimal subfield of $L$ containing $S$.
Let $k_0$ be the subfield $(L_0)^G$ fixed by $G$.
The action of $G$\/ on $S$ is faithful, so $G$\/ acts faithfully on $L_0$,
and $L_0/k_0$ is Galois with group $G$.
In this way we reduce to the case 
where $k$ and $L$ are finitely generated fields
(finitely generated over their minimal subfield).

Choose finitely generated $\Z$-algebras $A \subseteq B$
with fraction fields $k$ and $L$, respectively.
Without loss of generality  
we may assume, by localization, that $B$ is a finite \'etale Galois algebra over $A$.
Since $L$ is not algebraic over a finite field,
$\dim A = \dim B \geq 1$.
By~\cite[Theorem~4]{poonen-residue}, there is a maximal
ideal $\mm_1$ of $B$ lying over a maximal ideal $\mm$ of $A$
such that the residue field extension $B/\mm_1$ over $A/\mm$ is trivial.
Thus $\mm$ splits completely: 
if $n=\#G$,
there are $n$ distinct maximal ideals $\mm_1,\dots,\mm_n$ 
of $B$ lying over $\mm$,
and they are are permuted transitively by $G$.
By~\cite[Proposition~1.11]{atiyah-macdonald}, 
there exists a nonzero $\be \in \mm_1$
lying outside all of $\mm_2,\dots,\mm_n$.
We can label the conjugates $\be_i$ of~$\be$
so that $\be_i \in \mm_j$ if and only if $i=j$.
Any nontrivial relation $\prod_{i=1}^n \be_i^{b_i}= 1$ with $b_i \in \Z$,
would, after moving the factors with negative exponent to the other
side, give an equality between an element in $\mm_i$ 
and an element outside $\mm_i$, for some $i$.
Hence the ${\Z}G$-module generated by $\be$ in $L^*$ is free of rank~1.
\end{proof}

\begin{corollary}
\label{C:multiplicative S_r}
Let $k$ be a field that is not algebraic over a finite field.
If $k$ has a Galois extension with Galois group $S_r$,
then there exists $\al \in (k^s)^*$ with $r(\al)=d(\al)=r$.
\end{corollary}

\begin{proof}
Let $L$ be the $S_r$-extension of $k$.
By Proposition~\ref{P:multiplicative group contains ZG},
the $\Z{S_r}$-module $L^*$ contains a copy of $\Z{S_r}$,
which contains a copy of the $\Z{S_r}$-module $\Z^r$
on which $S_r$ acts by permuting coordinates.
The element $(1,0,\dots,0) \in \Z^r$ corresponds to $\al \in L^*$
with the desired properties.
\end{proof}

\begin{corollary}
\label{C:multiplicative representations}
Let $k$ be a field that is not algebraic over a finite field,
and let $G$ be a finite group.
Suppose that $G=\Gal(K/k)$ for some Galois extension $K$ of $k$,
and that there is a faithful $r$-dimensional subrepresentation $V$
of the regular representation of $G$ over~$\Q$.
Then there exists $\al \in K^*$ whose conjugates generate a torsion-free multiplicative group with $r(\al)=r$ and $d(\al)=[K:k]=\#G$.
\end{corollary}

\begin{proof}
Apply Proposition~\ref{P:multiplicative group contains ZG}
and then Lemma~\ref{L:faithful subrep} with $k=\Q$.
This gives $\al \in K^* \tensor_\Z \Q$ with the desired properties,
and we replace $\al$ by a power so that it is represented
by an element of $K^*$.
\end{proof}

We now prove the final statement of Theorem~\ref{T:mult}.
Since $k$ is Hilbertian, $k$ has $S_r$-extensions for all $r$.
In particular, $k$ is not algebraic over a finite field.
Applying Corollary~\ref{C:multiplicative S_r} yields
$\al$ with $r(\al)=d(\al)=r$.
Combining Corollaries \ref{C:inverseGalois} 
and~\ref{C:multiplicative representations}
gives a different $\al$ with 
$r(\al)=r$ and $d(\al)=\dmax(r)$, for any $r \geq 1$.
\end{proof}

We end by giving an explicit algebraic number of
conjugate rank~$n$ and degree $2^n n!$ over~$\Q$.

\begin{proposition}
Let $\sqrt{r_1}$, \dots, $\sqrt{r_n}$ be
as in Proposition~\ref{P:nonexceptionalexample}.
Let $s_i = (1+\sqrt{r_i})/(1-\sqrt{r_i})$ and
$\al=s_1 s_2^2 \cdots s_n^n$.
Then $r(\al)=n$ and $d(\al)=2^n n!$ over $\Q$.
\end{proposition}

\begin{proof}
The proof of Proposition~\ref{P:nonexceptionalexample}
showed that $[\Q(\sqrt{r_1},\dots,\sqrt{r_n}\,):\Q]=2^n n!$,
so its Galois group $G$\/ is the signed permutation group $W(B_n)$.
The elements of $G$\/ act on $\al$ by permuting the exponents $1,2,\dots,n$
and changing their signs independently.
In particular, the group generated by the conjugates of $\al$ 
is of finite index in the subgroup generated
by the $s_i$.
On the other hand, the $s_i$ are multiplicatively independent
since they are not roots of unity and
since there is an automorphism inverting any one of them
while fixing all the others.
Thus $\al$ has $2^n n!$ distinct conjugates,
and they generate a subgroup of rank $n$.
\end{proof}


\section*{Acknowledgments}

\noindent 
We thank Hendrik Lenstra for suggesting the proof of Lemma \ref{L:diagonal action}, and Walter Feit for some helpful correspondence.
N.D.E. also thanks David Moulton for communicating the topic to him, as a problem
on the list compiled by Gerry Myerson at the 2001 Asilomar meeting.
N.B. was supported by a UK EPSRC postgraduate award.
A.D.\ was supported in part by 
the Lithuanian State Studies and Science Foundation, 
and the London Mathematical Society.
N.D.E.\ was supported in part by NSF grant DMS-0200687.
B.P.\ was supported in part by NSF grant DMS-0301280 and a Packard Fellowship.


\bibliographystyle{amsalpha}
\bibliography{bjorn}

\end{document}